\documentclass[12pt]{article}
\usepackage{amsfonts,amssymb,version}

\begin{document}

\centerline{\large\bf Smooth homomorphisms admit noetherian reductions}

\bigskip

\centerline{Nitin Nitsure}

\begin{abstract}
  We give a short proof that any smooth (means formally smooth and finitely presented)
  homomorphism of rings can 
  be obtained by base change from a smooth homomorphism of noetherian rings.
  Together with the elegant short proof by J. Conde-Lago that smooth homomorphisms
  of noetherian rings are flat, this gives a short and elementary proof
  of the theorem of Grothendieck that all smooth homomorphisms are flat.
\end{abstract}

\centerline{MSC-2020: 14B25, 13B40.}

\medskip

A fundamental theorem of Grothendieck [2] 
says that any smooth (means formally smooth and finitely presented) homomorphism 
of commutative rings is flat. The proof given in [2] is 
quite intricate. Recently, J. Conde-Lago [1] gave a short and simple proof
that a smooth homomorphism of noetherian rings
is flat. The noetherian hypothesis is crucial to the proof of Conde-Lago, as it
enables the use of the standard machinery of completions of noetherian rings.

The main result of this note is the following, which allows us to reduce
arbitrary smooth homomorphisms to smooth homomorphisms between noetherian rings.

\medskip

{\bf Proposition 1.} {\it Let $A\to B$ be a ring homomorphism that is smooth, that is,
  finitely presented and formally smooth. Then there exist
finite-type $\mathbb Z$-algebras $A_0$ and $B_0$, a smooth homomorphism
$A_0 \to B_0$ and a homomorphism $A_0\to A$, such that 
$A\otimes_{A_0}B_0$ is isomorphic to $B$ as an $A$-algebra.

Moreover, we can choose $A_0 \to A$ to be injective, in which case 
the induced homomorphism $B_0 \to B$ is also injective.}

\medskip

{\bf Proof.} As $B$ is finitely presented over $A$, we can assume that
$B = P/I$ where $P$ is a polynomial ring $A[x_1,\ldots, x_n]$ in finitely
many variables, 
and $I\subset P$ is a finitely generated ideal $(f_1,\ldots, f_m)$.
Let $p: P \to P/I$, $p': P\to P/I^2$, and $q: P/I^2\to P/I$ denote
the quotient maps. 
Let $b_i = p(x_i) \in P/I = B$. 
As by assumption $B$ is formally smooth over $A$,  
the quotient homomorphism $q: P/I^2 \to P/I = B$ admits a section
$s: B \to P/I^2$ over $A$.
As $s$ is a section, there exist elements $g_1,\ldots, g_n \in I$ such that
$s(b_i) = p'(x_i + g_i) \in P/I^2$ for all $i$. 
As $g_i\in I$, there exist $u_{ij}\in P$ for $1\le i\le n$
and $1\le j\le m$ such that $g_i = \sum_j\, u_{ij}f_j$.
As $0 = p(f_j) = f_j(b_1,\ldots,b_n)  \in P/I$, applying $s$ we have
$0 = sp(f_j) = sf_j(b_1,\ldots,b_n) = 
f_j(s(b_1),\ldots,s(b_n)) =
f_j(p'(x_1+g_1),\ldots, p'(x_n+g_n))
= p'f_j(x_1+g_1,\ldots,x_n+g_n) \in P/I^2$, so  
$f_j(x_1+g_1,\ldots,x_n+g_n)\in {\mathop{\rm ker}\nolimits}(p') = I^2$. Hence there exist elements
$h_{jk\ell}\in P$ for $1\le j,k,\ell \le m$ such that 
$f_j(x_1+g_1,\ldots,x_n+g_n) = \sum_{k\ell} \, h_{jk\ell}f_kf_{\ell}$. 

Now let $A_0$ be the subring (sub-$\mathbb Z$-algebra) of $A$ generated
by all coefficients of the polynomials $f_j$, $u_{ij}$ and $h_{jk\ell}$,
with inclusion homomorphism $A_0 \hookrightarrow A$.
This ring is noetherian by the Hilbert basis theorem.
Let $P_0 = A_0[x_1,\ldots, x_n]\subset P$, which is a polynomial
ring over $A_0$ in the variables $x_i$'s,
as they remain algebraically independent over $A_0$.
Note that by definition of $A_0$, the elements 
$f_j$, $u_{ij}$ and $h_{jk\ell}$ of $P$ lie in $P_0$. 
The ideal $I_0 = \sum_j\, P_0f_j\subset P_0$ is contained in $I$, and 
$I = PI_0$. 
Let $B_0 = P_0/I_0$, which is a finitely presented $A_0$-algebra,
and let $B_0 \to B$ be the $A_0$-algebra homomorphism
induced by the inclusion $P_0\hookrightarrow P$.
Let $p_0: P_0 \to P_0/I_0$, $p'_0: P_0\to P_0/I_0^2$,
and $q_0: P_0/I_0^2\to P_0/I_0$ denote
the quotient maps. 
This induces
an $A_0$-algebra homomorphism $A\otimes_{A_0}B_0 \to B$, which as an isomorphism 
because the natural maps $A\otimes_{A_0}P_0 \to P$ and 
$P\otimes_{P_0}(P_0/I_0) \to P/I$ are isomorphisms. 

We now define a section $s_0 : B_0 = P_0/I_0 \to P_0/I_0^2$ of the
quotient $q_0: P_0/I_0^2 \to P_0/I_0 = B_0$ over $A_0$.
As $u_{ij}\in P_0$, we have $g_i = \sum_j\, u_{ij}f_j\in I_0$. We have a unique 
$A_0$-homomorphism $\sigma: P_0 \to P_0$ defined by
$\sigma(x_i) = x_i + g_i$. As $f_j,\,g_i,\, h_{jk\ell} \in P_0$, 
both sides of the equality 
$f_j(x_1+g_1,\ldots,x_n+g_n) = \sum_{k\ell} \, h_{jk\ell}f_kf_{\ell}$
are in $P_0$. As this equality is satisfied in $P$ it is
satisfied in $P_0\subset P$.
Hence $\sigma(f_j) = f_j(x_1+g_1,\ldots,x_n+g_n) \in I_0^2$. 
This shows that $\sigma: P\to P$ satisfies $\sigma(I_0)\subset I_0^2$, and so
$\sigma$ descends to an $A_0$-algebra homomorphism
$s_0 : B_0 = P_0/I_0 \to P_0/I_0^2$. 

By definition, $s_0p_0 (x_i) = p_0'(x_i + g_i)  \in P_0/I_0^2$, and
$q_0\circ p'_0 = p_0$.
Hence
$q_0s_0p_0(x_i) = q_0p_0'(x_i  + g_i) = p_0(x_i + g_i) = p_0(x_i)$ as $p_0(g_i) = 0$. 
This shows that $s_0 : B_0 \to P_0/I_0^2$ is a section of $q_0 : P_0/I_0^2 \to P_0/I_0 = B_0$.
As $P_0$ is formally smooth over $A_0$, this implies that $A_0 \to B_0$ is formally smooth,
and hence smooth, being finitely presented.

By our choice, $A_0$ is a subring of $A$, with $A_0\to A$ the inclusion. 
By the elementary proof by Conde-Lago [1],
any smooth homomorphism $A_0 \to B_0$ of noetherian rings 
is flat. The homomorphism $B_0 \to B$ is the base-change of $A_0 \to A$
under the flat homomorphism $A_0 \to B_0$, and so it is again injective if $A_0\to A$ is injective. 
\hfill$\square$

\bigskip

{\bf Corollary 2.} Any smooth (means formally smooth and finitely presented)
homomorphism of rings $A\to B$ is flat.

\medskip

{\bf Proof.} With notation as in the above proposition, we have shown that $A\to B$ is
the base change of a smooth homomorphism $A_0\to B_0$ of noetherian rings.  
By the elementary proof by Conde-Lago [1], $A_0\to B_0$ is flat.
Hence its base change $A\to B$ is flat.
\hfill$\square$

\bigskip

{\bf References}

[1] J. Conde-Lago: A short proof of smooth implies flat. 
Communications in Algebra 2017, vol 45, 774-775.

[2] A. Grothendieck and J. Dieudonn\'e: 
\'El\'ements de G\'eom\'etrie Alg\'ebrique IV. Publ. Math. IHES, vols. 20, 24, 28, 32 (1964-67).

\bigskip

{\small 
Nitin Nitsure \\
Retired Professor, School of Mathematics \\
Tata Institute of Fundamental Research \\
Mumbai 400 005, India \\
email: nitsure@gmail.com
}

\end{document}